\documentclass[amstex,12pt,amssymb]{article}

\usepackage{mathtext}
\usepackage[cp1251]{inputenc}
\usepackage[T2A]{fontenc}
\usepackage[dvips]{graphicx}
\usepackage{amsmath}
\usepackage{amssymb}
\usepackage{amsxtra}
\usepackage{latexsym}
\usepackage{ifthen}

\newcounter{lemma}[section]

\newcounter{corol}[section]

\newcounter{rem}[section]

\newcounter{theorem}[section]

\newcounter{proposition}[section]

\numberwithin{equation}{section}

\begin{document}

\markboth{\centerline{E. SEVOST'YANOV}}{\centerline{THE RADIUS OF
INJECTIVITY OF LOCAL RING $Q$--HOMEOMORPHISMS}}

\def\cc{\setcounter{equation}{0}
\setcounter{figure}{0}\setcounter{table}{0}}

\overfullrule=0pt

%\normalsize\large

\author{{E. SEVOST'YANOV}\\}

\title{
{\bf THE RADIUS OF INJECTIVITY OF LOCAL RING $Q$--HOMEOMORPHISMS}}

\date{\today}
\maketitle

%\large
\begin{abstract} The paper is devoted to the study
of mappings with non--bounded characteristics of quasiconformality.
The analog of the theorem about radius injectivity of locally
quasiconformal mappings was proved for some class of mappings. There
are found sharp conditions under which the so called local
$Q$--homeomorphisms are injective in some neighborhood of a fixed
point.
\end{abstract}


\section{Introduction}

Here are some definitions. Everywhere below, $D$ is a domain in
${\Bbb R}^n,$ $n\ge 2,$ $m$ be a measure of Lebesgue in ${\Bbb
R}^n,$ and ${\rm dist\,}(A,B)$ is the Euclidean distance between the
sets $A$ and $B$ in ${\Bbb R}^n.$ The notation $f:D\rightarrow {\Bbb
R}^n$ assumes that $f$ is continuous on its domain. In what follows
$(x,y)$ denotes the standard scalar multiplication of the vectors
$x,y\in {\Bbb R}^n,$ ${\rm diam\,}A$ is Euclidean diameter of the
set $A\subset {\Bbb R}^n,$
$$B(x_0, r)=\left\{x\in{\Bbb R}^n: |x-x_0|< r\right\}\,,\quad {\Bbb B}^n
:= B(0, 1)\,,$$ $$S(x_0,r) = \{ x\,\in\,{\Bbb R}^n :
|x-x_0|=r\}\,,\quad{\Bbb S}^{n-1}:=S(0, 1)\,,$$ $\omega_{n-1}$
denotes the square of the unit sphere ${\Bbb S}^{n-1}$ in ${\Bbb
R}^n,$ $\Omega_{n}$ is a volume of the unit ball ${\Bbb B}^{n}$ in
${\Bbb R}^n.$ A mapping $f:D\rightarrow {\Bbb R}^n$ is said to a
{\it local homeomorphism} if for every $x_0\in D$ there is a number
$\delta>0$ such that a mapping $f|_{B(x_0, \delta)}$ to be a
homeomorphism.

\medskip
Recall that a mapping $f:D\rightarrow {\Bbb R}^n$ is said to be {\it
a mapping with bounded distortion}, if the following conditions
hold:

\noindent 1) $f\in W_{loc}^{1,n},$

\noindent 2) a Jacobian $J(x,f):={\rm det\,}f^{\,\prime}(x) $ of the
mapping $f$ at the point $x\in D$ preserves the sign almost
everywhere in $D,$

\noindent 3) $\Vert f^{\,\prime}(x) \Vert^n \le K \cdot |J(x,f)|$ at
a.e. $x\in D$ and some constant $K<\infty,$ where $$\Vert
f^{\,\prime}(x)\Vert:=\sup\limits_{h\in {\Bbb R}^n:
|h|=1}|f^{\,\prime}(x)h|\,,$$
see., e.g., \cite[\S\,3, Ch.~I ]{Re$_2$}, or definition 2.1 of the
section 2 Ch. I in \cite{Ri}. In this case we also say that $f$ is
{\it $K$--qua\-si\-re\-gu\-lar,} where $K$ is from condition 3)
meaning above.

\medskip
The following result was proved in the work \cite{MRV$_3$} by O.
Martio, S. Rickman and J.~V\"{a}is\"{a}l\"{a}, see
\cite[Theorem~2.3]{MRV$_3$} or \cite[Theorem~3.4.III]{Ri}, see also
the paper \cite{KOR}.

\medskip
{\bf Statement 1.} {\sl If $n\ge 3$ and $f:{\Bbb B}^n\rightarrow
{\Bbb R}^n$ is a $K$--qua\-si\-re\-gu\-lar local homeomorphism, then
$f$ is injective in a ball $B\left(0, \psi(n, K)\right),$ where
$\psi$ is a positive number depending only on $n$ and $K.$}

\medskip
A goal of the present paper is a proof of the analog of the
Statement 1 for more general classes of mappings of ring
$Q$--ho\-me\-o\-mor\-phisms. To introduce this class of the
mappings, we give some definitions.

\medskip
A curve $\gamma$ in ${\Bbb R}^n$ is a continuous mapping $\gamma
:\Delta\rightarrow{\Bbb R}^n$ where $\Delta$ is an open, closed or
half--open interval in ${\Bbb R}.$ Given a family $\Gamma$ of paths
$\gamma$ in ${\Bbb R}^n,$ $n\ge 2,$ a Borel function $\rho:{\Bbb
R}^n \rightarrow [0,\infty]$ is called {\it admissible} for
$\Gamma,$ abbr. $\rho \in {\rm adm}\,\Gamma,$ if
$$\int\limits_{\gamma} \rho(x)\,|dx| \ge 1$$ for each
$\gamma\in\Gamma .$ The {\it modulus} of $\Gamma$ is the quantity
$$M(\Gamma) =\inf\limits_{ \rho \in {\rm adm}\,\Gamma}
\int\limits_{{\Bbb R}^n} \rho^n(x)\, dm(x)\,.$$

\medskip
Given a domain $D$ and two sets $E$ and $F$ in $\overline{{\Bbb
R}^n},$ $n\ge 2,$ $\Gamma(E,F,D)$ denotes the family of all paths
$\gamma:[a,b] \rightarrow {\overline{{\Bbb R}^n}}$ which join $E$
and $F$ in $D,$ i.e., $\gamma(a)\in E,$ $\gamma(b)\in F$ and
$\gamma(t) \in D$ for $a<t<b.$

Let $D$ be a domain in ${{\Bbb R}^n},$  $Q: D \rightarrow
[0,\infty]$ be a (Lebesgue) measurable function. Set
$$A(x_0,r_1,r_2) = \{ x\in{\Bbb R}^n : r_1<|x-x_0|<r_2\}\,.$$
We say that a mapping $\overline{{\Bbb R}^n}$ is a {\it ring
$Q$--mapping at a point} $x_0\in D$ if
\begin{equation}\label{1.8}
M\left(f(\Gamma(S(x_0, r_1), S(x_0, r_2), A(x_0,r_1,r_2)))\right)
\le \int\limits_{A(x_0,r_1,r_2)} Q(x)\cdot \eta^n(|x-x_0|)\ dm(x)
\end{equation}
for every ring $A(x_0,r_1,r_2),$ $0<r_1<r_2< r_0={\rm
dist}(x_0,\partial D),$ and for every Lebesgue measurable function
$\eta: (r_1,r_2)\rightarrow [0,\infty ]$ such that
$$\int\limits_{r_1}^{r_2}\eta(r)dr\ge 1\,.$$
If the condition (\ref{1.8}) holds at every point $x_0\in D,$ then
we also say that $f$ is a ring $Q$--mapping in the domain $D,$ see
\cite[section~7]{MRSY$_2$}.

\medskip In what follows
$q_{x_0}(r)$ denotes the integral average of $Q(x)$ under the sphere
$|x-x_0|=r,$
\begin{equation}\label{eq17}
q_{x_0}(r):=\frac{1}{\omega_{n-1}r^{n-1}}\int\limits_{|x-x_0|=r}Q(x)\,dS\,,
\end{equation}
where $dS$ is element of the square of the surface $S.$

\medskip
One of the main results of the paper is following.

\medskip
\begin{theorem}\label{th1A} {\sl Let $n\ge 3$ and $f:{\Bbb B}^n\rightarrow {\Bbb R}^n$ is a
local ring $Q$--ho\-me\-o\-mor\-phism at the point $x_0=0,$ such
that $Q\in L_{loc}^{1}({\Bbb B}^n)$ and
\begin{equation}\label{eq1A}
\int\limits_{0}^{1}\frac{dt}{tq_0^{1/(n-1)}(t)}=\infty\,.
\end{equation}
Then $f$ is injective in a ball $B\left(0, \delta(n, Q)\right),$
where $\delta$ is a positive number depending only on $n$ and
function $Q.$ From other hand, the condition (\ref{eq1A}) is
precise, in fact, for every $\delta>0$ and every $Q\in
L_{loc}^{1}({\Bbb B}^n)$ with $Q(x)\ge 1$ a.e. and
\begin{equation}\label{eq1B}
\int\limits_{0}^{1}\frac{dt}{tq_0^{1/(n-1)}(t)}<\infty
\end{equation}
there exists a mapping $f=f_Q:{\Bbb B}^n\rightarrow {\Bbb R}^n$
which is local ring $Q$--homeomorphism at the point $x_0=0$ and
which is not injective in $B(0, \delta).$}
\end{theorem}

%For the proof of the Theorem meaning above we have used a part of
%the sketch of the proof used in \cite{MRV$_3$}, see also \cite{Ri},
%and \cite{KOR}.

\section{The main Lemma}
\setcounter{equation}{0}

\medskip
A set $Q\subset\overline{{\Bbb R}^n}$ is said to be {\it relatively
locally connected} if every point in $\overline{Q}$ has arbitrary
small neighborhoods $U$ such that $U\cap Q$ is connected.

We also need following statements, see \cite[Lemmas~3.1.III--
3.3.III]{Ri}.

\medskip
\begin{proposition}\label{pr1A}{\sl\,
Let $f:G\rightarrow\overline{{\Bbb R}^n}$ be a local homeomorphism,
let $Q$ be a simply connected and locally pathwise connected set in
$\overline{{\Bbb R}^n},$ and let $P$ be a component of $f^{\,-1}(Q)$
such that $\overline{P}\subset G.$ Then $f$ maps $P$
homeomorphically onto $Q.$ If, in addition, $Q$ is relatively
locally connected, $f$ maps $\overline{P}$ homeomorphically onto
$\overline{Q}.$}
\end{proposition}

\medskip
\begin{proposition}\label{pr2A}{\sl\,
Let $f:G\rightarrow\overline{{\Bbb R}^n}$ be a local homeomorphism
and let $F$ be a compact set in $G$ such that $f_F$ is injective.
Then $f$ is injective in a neighborhood of $F.$}
\end{proposition}

\medskip
\begin{proposition}\label{pr3A}{\sl\,
Let $f:G\rightarrow\overline{{\Bbb R}^n}$ be a local homeomorphism,
let $A, B\subset G,$ and let $f$ be homeomorphic in $A$ and $B.$ If
$A\cap B\ne \varnothing$ and $f(A)\cap f(B)$ is connected, then $f$
is homeomorphic in $A\cup B.$ }
\end{proposition}

\medskip
Finally, we need the following statement of P. Koskela, J, Onninen
and K. Rajala, see \cite[Lemma~3.1]{KOR}.

\medskip
\begin{proposition}\label{pr4}
{\sl Let $n\ge 1$ and $r>0.$ Let $a\ne b,$ $a, b\in S(0, r).$ Then
there exists a point $p=p(a, b)\in B(0, r)$ such that for every
$t\in \left(\frac{r}{2}, \frac{\sqrt{3}r}{2}\right)$ either
$$0, b\in B(p, t)\quad and \quad a\not\in B(p, t)$$
or
$$a, b\in B(p, t)\quad and \quad 0\not\in B(p, t)\,.$$ }
\end{proposition}

\medskip
The following Lemma plays the main role in the following.

\medskip
\begin{lemma}\label{lem2}{\sl\, Let $n\ge 3,$ $Q:{\Bbb B}^n\rightarrow
[0, \infty]$ and $f:{\Bbb B}^n\rightarrow {\Bbb R}^n$ is a local
ring $Q$--ho\-me\-o\-mor\-phism at the point $x_0=0.$ Suppose that
there exist a function $\psi:(0,1)\rightarrow [0, \infty]$ and a
constant $C=C(n, Q, \psi)$ such that
\begin{equation}\label{eq2A}
0<I(r_1, r_2):=\int\limits_{r_1}^{r_2}\psi(t)dt<\infty \quad\forall
\,r_1, r_2\in (0,1)\end{equation}
and for some $\alpha>0$
\begin{equation}\label{eq3A}
\int\limits_{r_1<|x|<r_2}Q(x)\psi(|x|)dm(x)\le C\cdot
I^{n-\alpha}(r_1, r_2)\,.
\end{equation}
Let \begin{equation}\label{eq4A} I(0,
1):=\int\limits_{0}^{1}\psi(t)dt=\infty\,,\end{equation} then $f$ is
injective in a ball $B\left(0, \delta(n, Q, \psi)\right),$ where
$\delta$ is a positive number depending only on $n,$ functions $Q$
and $\psi.$ }
\end{lemma}

\medskip
\begin{proof}
{\bf The 1 step.} We may assume $f(0)=0.$ Let $r_0=\sup\{r\in {\Bbb
R}: r>0, \overline{U(0, r)}\subset {\Bbb B}^n\},$ where $U(0, r)$ is
the $0$--component of $f^{\,-1}(B(0, r)).$ Clearly $r_0>0.$
Fix $r<r_0$ and set $U=U(0, r),$
$$l^{\,*}=l^{\,*}(0, f, r)=\inf\{|z|: z\in \partial U\}\,,$$
$$L^{\,*}=L^{\,*}(0, f, r)=\sup\{|z|: z\in \partial U\}\,.$$
By Proposition \ref{pr1A}, $f$ maps $\overline{U}$ homeomorphically
onto $\overline{B(0, r)}.$ Thus $f$ is injective in $B(0, l^{\,*})$
and it suffices to find a lower bound for $l^{\,*}.$

\medskip
{\bf The 2 step.} Note that $L^{\,*}\rightarrow 1$ as $r\rightarrow
r_0.$ Suppose the contrary: $L^{\,*}\not\rightarrow 1$ as
$r\rightarrow r_0.$

\medskip
а) Remark that  $U(0, r_1)\subset U(0, r_2)$ as $0<r_1<r_2<r_0.$ In
fact, let us assume that there exists $x\in U(0, r_1)\setminus U(0,
r_2).$ Since $f(U(0, r_i))=B(0, r_i),$ $i=1,2,$ we have $f(x)=y\in
B(0, r_1)$ and $f(z)=y\in B(0, r_1),$ $z\ne x.$ However, this
contradicts to the Proposition \ref{pr3A}, because $f$ is
homeomorphism in $U(0, r_1)\cup U(0, r_2)$ in this case.

b) It follows from a) that the function $L^{\,*}$ is increase by $r$
and, consequently, there exists the limit of $L^{\,*}$ as
$r\rightarrow r_0.$ Then $L^{\,*}\rightarrow \varepsilon_0$ as
$r\rightarrow r_0,$ where $\varepsilon_0\in (0, 1).$ In this case,
$U(0, r)\subset B(0, \varepsilon_0)$ for every $0<r<r_0.$

c) Remark that $B(0, r_0)\subset f(B(0, \varepsilon_0)).$ In fact,
let $y\in B(0, r_0),$ then $y\in B(0, r_1)$ for some $r_1\in (0,
r_0).$ It follows from hence that there exists $x\in U(0, r_1)$ with
$f(x)=y$ and, consequently, $y\in f(B(0, \varepsilon_0)),$ i.e.,
$B(0, r_0)\subset f(B(0, \varepsilon_0)).$

d) Remark that $\overline{B(0, r_0)}\subset f(\overline{B(0,
\varepsilon_0)})$ and, consequently, by the openness of $f,$ $f(B(0,
\varepsilon_1))$ contains some neighborhood of $\overline{B(0,
r_0)}$ for every $\varepsilon_1\in (0, \varepsilon_0).$ Thus, $U(0,
r_0)$ lies inside of $\overline{B(0, \varepsilon_0)},$ that
contradicts to the definition of $r_0.$ The contradiction obtained
above implies that $L^{\,*}\rightarrow 1$ as $r\rightarrow r_0$ that
is desired conclusion.

\medskip
{\bf The 3 step.} Pick $x$ and $y\in \partial U$ such that
$|x|=L^{\,*}$ and $|y|=l^{\,*}.$ Note that, by the definition of
$U,$ $f(x),$ $f(y)\in S(0, r).$ By Proposition \ref{pr4} there
exists a point $p\in B(0, r)$ such that, for every $t\in
\left(\frac{r}{2}, \frac{\sqrt{3}r}{2}\right),$ $f(x)\in B(p, t)$
and either $0\in B(p, t)$ and $f(y)\not\in B(p, t),$ or $0\not\in
B(p, t)$ and $f(y)\in B(p, t).$ Fix such a $t.$ Note that $0, f(y)$
and $f(x)\in \overline{f(B(0, l^{\,*}))}$ and, consequently, $f(B(0,
l^{\,*}))\cap B(p, t)\ne\varnothing\ne f(B(0, l^{\,*}))\setminus
B(p, t).$ Since $f(B(0, l^{\,*}))$ is connected, this implies that
there exists a point $z_t\in S(p, t)\cap  f(B(0, l^{\,*})),$ see
\cite[Theorem~1.I.46.5]{Ku$_2$}.

\medskip
Let $z_t^{\,*}$ be the unique point in $f^{\,-1}(z_t)\cap B(0,
l^{\,*}).$ Let $C_t(\varphi)\subset S(p, t)$ be the spherical cap
with center $z_t$ and opening angle $\varphi,$
$$C_t(\varphi)=\{y\in {\Bbb R}^n: |y-p|=t, (z_t-p, y-p)>t^{\,2}\cos\varphi\}\,.$$
Let $\varphi_t$ be the supremum of all $\varphi$ for which the
$z_t^{\,*}$--component of $f^{\,-1}(C_t(\varphi))$ gets mapped
homeomorphically onto $C_t(\varphi).$ Let $C_t=C_t(\varphi_t)$ and
let $C_t^{\,*}$ be the $z_t^{\,*}$--component of $f^{\,-1}(C_t).$

\medskip
{\bf The 4 step.} We claim that $C_t^{\,*}$ meets $S(0, L^{\,*}).$
Suppose this is not true.
\medskip

a) Since $C_t^{\,*}$ is connected and $C_t^{\,*}\cap B(0,
L^{\,*})\ne \varnothing,$ this implies that $C_t^{\,*}\subset B(0,
L^{\,*}),$ see \cite[Theorem~1.I.46.5]{Ku$_2$}. Remark that, in this
case, $C_t^{\,*}$ is a compact subset of $U$ and by Proposition
\ref{pr1A} $f$ maps $\overline{C_t^{\,*}}$ homeomorphically onto
$\overline{C_t}.$ (It is not true at $n=2$ because $C_t(\pi)$ is not
relatively locally connected). By Proposition \ref{pr2A} $f$ is
injective in a neighborhood of $\overline{C_t^{\,*}}.$ Thus
$\varphi_t=\pi,$ $\overline{C_t}=S(p, t)$ and $\overline{C_t^{\,*}}$
is a topological $(n-1)$--sphere in ${\Bbb R}^n.$ Note that bounded
component $D$ of ${\Bbb R}^n\setminus \overline{C_t^{\,*}}$
contained in $B(0, L^{\,*}).$ Now $\overline{f(D)}$ is a compact
subset of $f({\Bbb B}^n)$ and, since the mapping $f$ is open,
$\partial f(D)\subset f(\partial D)=S(p, t).$

\medskip
b) Remark that $f(D)\subset B(p, t).$ In fact, let $f(D)\not\subset
B(p, t),$ then there exists $y\in f(D)\setminus \overline{B(p ,t)}.$
Now we have $\left(f({\Bbb B}^n)\setminus \overline{B(p,
t)}\right)\cap f(D)\ne\varnothing$ and, since $f(D)$ is compact
subdomain of $f({\Bbb B}^n),$ $\left(f({\Bbb B}^n)\setminus
\overline{B(p, t)}\right)\setminus f(D)\ne\varnothing.$ Since
$f({\Bbb B}^n)\setminus \overline{B(p, t)}$ is connected, this
implies that there exists $z\in \partial f(D)\cap \left(f({\Bbb
B}^n)\setminus \overline{B(p, t)}\right),$ see \cite[Theorem
~1.I.46.5]{Ku$_2$}, that contradicts to the inclusion $\partial
f(D)\subset S(p, t).$

c) Now $f(D)\subset B(p, t).$ Remark that $B(p, t)\subset f(D).$
Indeed, let there exists $a\in B(p, t)\setminus f(D).$ Since $B(p,
t)$ is connected and $B(p, t)\cap f(D)\ne\varnothing$ this implies
that $\partial f(D)\cap B(p, t)\ne \varnothing,$ see
\cite[Theorem~1.I.46.5]{Ku$_2$}. The last relation contradicts to
the inclusion $\partial f(D)\subset S(p, t).$

d) Thus $f(D)=B(p, t).$ By the definition, $D$ is a component of
$f^{\,-1}(B(p, t)).$ By Proposition \ref{pr1A} $f$ maps
$\overline{D}$ onto $\overline{B(p, t)}$ homeomorphically.

\medskip
e) Since $z_t^{\,*}\in \overline{C_t^{\,*}}\cap U,$
$\overline{D}\cap \overline{U}\ne\varnothing.$ Since $f$ maps
$\overline{U}$ homeomorphically onto $\overline{B(0, r)},$ $f$ is
injective in $\overline{U}\cup \overline{D}$ by Proposition
\ref{pr3A}. This is impossible, because in view of the equality
$f(D)=B(p, t)$ and that $f(x)\in B(p, t)$ there exists a point
$x_1\ne x,$ $x_1\in D,$ such that $f(x_1)=f(x).$ Thus $C_t^{\,*}$
meets $S(0, L^{\,*}).$

\medskip
{\bf The 5 step.} Let $k_t^{\,*}\in C_t^{\,*}\cap S(0, L^{\,*})$ and
$k_t=f(k_t^{\,*}).$ Let $\Gamma_t^{\,\prime}$ be the family of all
curves connecting $k_t$ and $z_t$ in $C_t.$ Moreover, let
$\Gamma^{\,\prime}$ be the union of the curve families
$\Gamma^{\,\prime}_t,$ $t\in \left(\frac{r}{2},
\frac{\sqrt{3}r}{2}\right).$ Denote by $f_t$ the restriction of $f$
to $C_t^{\,*}.$ Then $f_t$ maps $C_t^{\,*}$ homeomorphically onto
$C_t.$ Furthermore, denote
$$\Gamma=\bigcup\limits_{t\in \left(\frac{r}{2}, \frac{\sqrt{3}r}{2}\right)}
\left\{f_t^{\,-1}\circ \gamma: \gamma\in
\Gamma_t^{\,\prime}\right\}\,.$$
Since for every $t\in \left(\frac{r}{2},
\frac{\sqrt{3}r}{2}\right),$ $z_t^{\,*}\in B(0, l^{\,*})$ and
$k_t\in S(0, L^{\,*}),$ be the definition of ring $Q$--map\-ping we
have
\begin{equation}\label{eq5B}
M(f(\Gamma(S(0, l^{\,*}), S(0, L^{\,*}), A(0, l^{\,*}, L^{\,*} ))))
\le \int \limits_{A(0, l^{\,*}, L^{\,*}
)}Q(x)\cdot\eta^{\,n}(|x|)dm(x)
\end{equation}
for every function $\eta: (l^{\,*}, L^{\,*})\rightarrow [0,\infty ]$
with
$$\int\limits_{l^{\,*}}^{L^{\,*}}\eta(r)dr\ge 1\,.$$
Setting $\eta(t)=\psi(t)/I(l^{\,*}, L^{\,*}),$ where $\psi$ is the
function from the condition of Lemma, we observe that $\eta$
satisfies the above condition. Now from (\ref{eq3A}) and
(\ref{eq5B}) we obtain that
$$M(\Gamma^{\,\prime})=M(f(\Gamma(S(0, l^{\,*}), S(0, L^{\,*}), A(0,
l^{\,*}, L^{\,*} )))) \le$$
\begin{equation}\label{eq6B}
\le \int \limits_{A(0, \, l^{\,*}, \, L^{\,*}
)}Q(x)\cdot\eta^{\,n}(|x|)dm(x)\le C/I^{\,\alpha}(l^{\,*},
L^{\,*})\,.
\end{equation}
On other hand, by \cite[Theorem~10.2]{Va$_1$},
\begin{equation}\label{eq7B}
\int\limits_{S(p, t)}\rho^n(x)dS\ge \frac{C_n}{t}
\end{equation}
for every $\rho$ for which $\int\limits_{\gamma}\rho(x)|dx|\ge 1$
for every $\gamma\in \Gamma_t^{\,\prime}.$ The integration of
(\ref{eq7B}) over $t$ yields
\begin{equation}\label{eq8B}
M(\Gamma^{\,\prime})\ge C_n^{\,\prime}
\end{equation}
for some constant $C_n^{\,\prime}>0.$ We obtain from (\ref{eq6B})
and (\ref{eq8B}) that
\begin{equation}\label{eq9B}
C_n\le C/I^{\,\alpha}(l^{\,*}, L^{\,*})\le C/I^{\,\alpha}(l^{\,*}(0,
f, r_0), L^{\,*}(0, f, r))
\end{equation}
because $I(\varepsilon_1, \varepsilon_2)> I(\varepsilon_3,
\varepsilon_2)$ as $\varepsilon_3>\varepsilon_1.$ Letting into the
limit as $r\rightarrow r_0$ in (\ref{eq9B}), we have
\begin{equation}\label{eq10B}
C_n\le C/I^{\,\alpha}(l^{\,*}(0, f, r_0), 1)\,.
\end{equation}
First of all, from the (\ref{eq10B}) follows that $I(\varepsilon,
1)<\infty$ for every $\varepsilon\in (0, 1).$ Follow, let
$l^{\,*}(0, f, r_0)\rightarrow 0,$ then it follows from (\ref{eq4A})
that the right hand of (\ref{eq10B}) tends to zero, that contradicts
to (\ref{eq10B}). Thus, $l^{\,*}(0, f, r_0)\ge \delta$ for all such
$f.$ The proof is complete.
\end{proof} $\Box$

\section{Proof of the main result}
\setcounter{equation}{0}

The following statement would be very useful, see \cite[Theorem
~3.15]{RS}.

\medskip
\begin{proposition}\label{pr1B}
{\sl Let $D$ be a domain in ${\Bbb R}^n,$ $n\ge 2,$ and $Q:D
\rightarrow [0,\infty ]$ a locally integrable measurable function. A
homeomorphism $f:D\rightarrow \overline{{\Bbb R}^n}$ is a ring
$Q$--homeomorphism at a point $x_0$ if and only if for every
$0<r_1<r_2< r_0= {\rm dist}\, (x_0,\partial D),$
$$M\left(\Gamma\left(f(S_1), f(S_2), f(D)\right)\right)\le
\frac{\omega_{n-1}}{I^{n-1}}$$
where $\omega_{n-1}$ is the area of the unit sphere in ${\Bbb R}^n,$
$q_{x_0}(r)$ is the average of $Q(x)$ over the sphere $|x-x_0|=r,$
$S_j=\{x\in{\Bbb R}^n: |x-x_0|=r_j\},$ $j=1,2,$ and
$$I=I(r_1,r_2)=\int\limits_{r_1}^{r_2}\
\frac{dr}{rq_{x_0}^{\frac{1}{n-1}}(r)}\,.$$}
\end{proposition}

\medskip
{\it Proof of Theorem \ref{th1A}}. Given $0<r_1<r_2<r_0=1$ consider
the function
\begin{equation}\label{eq1*****}\psi(t)\quad=\quad \left \{\begin{array}{rr}
1/[tq^{\frac{1}{n-1}}_{0}(t)]\ , & \ t\in (r_1, r_2)\ ,
\\ 0\ ,  &  \ t\notin (r_1, r_2)\ .
\end{array} \right.
\end{equation}
Note that $\psi$ satisfies all the conditions of Lemma \ref{lem2},
in particular,
$\int\limits_{r_1}^{r_2}\frac{dt}{tq_0^{1/(n-1)}(t)}<\infty$ by
\cite[Theorem~1]{Sev}, and by Fubini theorem,
$\int\limits_{r_1<|x|<r_2} Q(x)\cdot\psi^n(|x|)\ dm(x)\ =\
\omega_{n-1}\cdot I(r_1, r_2).$ Now the first part of the Theorem
follows from Lemma \ref{lem2}.

\medskip
To prove second part of the Theorem, we take $\delta>0$ and some
function $Q\in L_{loc}^{1}({\Bbb B}^n)$ satisfying (\ref{eq1B}). Set
$$f(x)=\frac{x}{|x|}\rho(|x|)\,,$$
where
$$\rho(r)=\exp\left\{-\int\limits_{r}^1\frac{dt}{t\widetilde{q}_{0}^{1/(n-1)}(t)}\right\}\,,
\qquad
\widetilde{q}_{0}(r):=\frac{1}{\omega_{n-1}r^{n-1}}\int\limits_{|x|=r}\widetilde{Q}(x)\,dS\,,
$$
$$\widetilde{Q}(x)\quad=\quad \left \{\begin{array}{rr} Q(x) , & \
|x|> \delta\ ,
\\ 1/K\ ,  &  |x|\le\delta\,,
\end{array} \right.$$
where $K\ge 1$ would be chosen bellow. Note that a mapping $f$ is a
ring $\widetilde{Q}$--homeomorphism at $x_0=0.$ In fact, we have
$f(S(0, r))=S(0, R),$ where
$R:=\exp\left\{-\int\limits_{r}^1\frac{dt}{t\widetilde{q}_{0}^{1/(n-1)}(t)}\right\}.$
Now $$f(\Gamma(S(0, r_1), S(0, r_2), A(0, r_1, r_2)))= \Gamma(S(0,
R_1), S(0, R_2), A(0, R_1, R_2))\,,$$ where
$R_i:=\exp\left\{-\int\limits_{r_i}^1\frac{dt}{t\widetilde{q}_{0}^{1/(n-1)}(t)}\right\},$
$i=1, 2.$ Now by \cite[section~7.5]{Va$_1$},
$$M(f(\Gamma(S(0, r_1), S(0, r_2), A(0, r_1, r_2))))=\frac{\omega_{n-1}}
{\left(\int\limits_{r_1}^{r_2}\frac{dt}{t\widetilde{q}_{0}^{1/(n-1)}(t)}\right)^{n-1}}\,.$$
Now $f$ is a ring $\widetilde{Q}$--homeomorphism at the point
$x_0=0$ by Proposition \ref{pr1B} and, consequently, is a ring
$Q$--mapping at $0.$ Note that under $\delta\rightarrow 0$ the image
$f(B(0, \delta))$ includes the ball $B(0, \sigma),$ where $\sigma$
does not depend on $\delta.$ Now we map the ball $B(0, \sigma)$ by
some map $g,$ which is $K$--quasiregular and local homeomorphism for
some $K\ge 1,$ but not injective in $B(0, \sigma)$; for instance,
let $g$ is a winding map, whose axes of rotation does not contain a
ball ${\Bbb B}^n=f(B(0, 1)),$  see \cite[section~5.1.I]{Re$_2$}.
Remark that $K$ does not depend on $\delta.$ Now we construct a
local ring $K\cdot Q(x)$--homeomorphism $f_2$ at zero, $f_2=g\circ
f,$ which is not injective in $B(0, \delta).$ Since $Q$ is arbitrary
locally integrable function with $Q\ge 1$ satisfying (\ref{eq1B}),
we can replace $Q$ on the $Q/K$ in the start of the second part of
the proof. So, we obtain a local ring $Q(x)$--homeomorphism with the
properties meaning above. The proof is complete. $\Box$

\section{Corollaries}

The following statement is a simple consequence from the first part
of the Theorem \ref{th1A}.

\medskip
\begin{corol}\label{cor1} {\sl Let $f:{\Bbb B}^n\rightarrow {\Bbb R}^n,$ $n\ge 3,$ be a local ring
$Q$--mapping at $x_0=0$ such that
\begin{equation}\label{eq11}
q_{0}(r)\le C\cdot\log^{n-1}\frac{1}{r}
\end{equation}
for some $C>0$ and $r\rightarrow 0.$
Then $f$ is injective in some ball $B\left(0, \delta(n, Q)\right)$
where $\delta$ depends only on $n$ and $Q.$ }
\end{corol}

\medskip
\begin{proof}
The desired conclusion follows from the Theorem \ref{th1A} in view
of (\ref{eq11}). In fact, it follows from the Fubini Theorem (see
\cite[Theorem~8.1, Ch.~III]{Sa}) that $Q\in L_{loc}^1({\Bbb B}^n),$
besides of that,  if follows from (\ref{eq11}) that (\ref{eq1A})
holds.
\end{proof} $\Box$

Following \cite{IR}, we say that a function $\varphi:D\rightarrow
{\Bbb R} $ has {\it finite mean oscillation} at a point $x_0 \in {D}
$ if
$$\limsup\limits_{\varepsilon\rightarrow 0}\frac{1}{\Omega_n\cdot \varepsilon^n}
\int\limits_{B(x_0 ,\varepsilon)}
|\varphi(x)-\widetilde{\varphi_{\varepsilon}}|dm(x)< \infty
$$
where
$$\widetilde{\varphi_{\varepsilon}}=
\frac{1}{\Omega_n\cdot \varepsilon^n}\int\limits_{B( x_0
,\varepsilon)} \varphi(x)\, dm(x)$$
%\end{equation}
%
is the average of the function $\varphi(x)$ over the ball
$B(x_0,\varepsilon)=\{x\in {\Bbb R}^n: |x-x_0|<\varepsilon\}.$

\medskip
We also say that a function $\varphi:D\rightarrow{\Bbb R}$ is of
finite mean oscillation in the domain $D,$ abbr. $\varphi\in FMO(D)$
or simply $\varphi\in FMO,$ if $\varphi$ has finite mean oscillation
at every point $x_0\in D.$ Note that $FMO$ is not $BMO_{loc},$ see
examples in \cite[p.~211]{MRSY$_2$}. It is well--known that
$L^{\,\infty }(D)\subset BMO(D)\subset L^p_{loc}(D)$ for all $1\le
p<\infty,$ see e.g. \cite{JN}, but $FMO(D)\not\subset L_{loc}^p(D)$
for any $p>1.$ The following statement can be found in
\cite[Lemma~6.1]{MRSY$_2$}.

\medskip
\begin{proposition}\label{pr7}{\sl\,
Let $0\in D\subset {\Bbb R}^n,$ $n\ge 3,$ $\varphi: D\rightarrow
{\Bbb R}$ be a nonnegative function having a finite mean oscillation
at $x_0=0.$ Then there exists $\varepsilon_0>0$ with
$$\int\limits_{B(0, \varepsilon_0)}\frac{\varphi(x)\,
dm(x)} {\left(|x| \log \frac{1}{|x|}\right)^n} <\infty\,.$$
%\end{equation}
%
%
}
\end{proposition}

\medskip
The following statement take a place.

\medskip
\begin{theorem}\label{th2} {\sl Let $g:{\Bbb B}^n\rightarrow {\Bbb R}^n,$ $n\ge 3,$ be a local
ring $Q$--mapping at $x_0=0$ such that $Q\in FMO(0).$ Then $g$ is
injective in some ball $B\left(0, \delta(n, Q)\right),$ where
$\delta$ is positive number depending only on $n$ and $Q.$ }
\end{theorem}

\medskip
\begin{proof}
Let $\varepsilon_0>0$ be a number from the Proposition \ref{pr7}.
Consider the mapping $f:=g(x\varepsilon_0),$ $x\in {\Bbb B}^n.$
Remark that $g$ be a local ring $Q(\varepsilon_0 x)$--mapping at
zero. Let us apply the Lemma \ref{lem2} for the mapping $g$ and
function $\psi=\frac{1}{\varepsilon_0t\log\frac{1}{\varepsilon_0
t}}.$ By Proposition \ref{pr7} we obtain that the relation
(\ref{eq3A}) holds at $\alpha=n$ for the function $\psi$ mentioned
above. Remark that, the relations (\ref{eq2A}) are (\ref{eq4A})
hold, also. The desired conclusion follows from the Lemma
\ref{lem2}. $\Box$
\end{proof}

\section{On normality of the families of local homeomorphisms}
\setcounter{equation}{0}

In what follows, we use in $\overline{{{\Bbb R}}^n}={{\Bbb
R}}^n\bigcup\{\infty\}$ the {\bf spherical (chordal) metric}
$h(x,y)=|\pi(x)-\pi(y)|$ where $\pi$ is the stereographic projection
of $\overline{{{\Bbb R}}^n}$ onto the sphere
$S^n(\frac{1}{2}e_{n+1},\frac{1}{2})$ in ${{\Bbb R}}^{n+1},$ i.e.
%\begin{equation}\label{eq6.2.15}
$$h(x,y)=\frac{|x-y|}{\sqrt{1+{|x|}^2} \sqrt{1+{|y|}^2}},\ \,\, x\ne
\infty\ne y, $$%\end{equation}
$$
h(x,\infty)=\frac{1}{\sqrt{1+{|x|}^2}}\ .
$$
It is clear that $\overline{{\Bbb R}^n}$ is homeomorphic to the unit
sphere ${\Bbb S}^n$ in ${\Bbb R}^{n+1}.$

The {\it spherical (chordal) diameter} of a set $E \subset
\overline{{\Bbb R}^n}$ is
%\begin{equation}\label{eq6.2.16}
$$h(E)=\sup_{x,y \in E} h(x,y)\,.$$

\medskip
Let $(X,d)$ and $\left(X^{{\,\prime}},{d}^{{\,\prime}}\right)$ be
metric spaces with distances $d$ and $d^{{\,\prime}}$, respectively.
A family $\frak{F}$ of continuous mappings $f:X\rightarrow
{X}^{{\,\prime}}$ is said to be {\it normal} if every sequence of
mappings $f_m \in \frak{F}$ has a subsequence $f_{m_k}$ converging
uniformly on each compact set $C \subset X$ to a continuous mapping.
Normality is closely related to the  following. A family $\frak{F}$
of mappings $f:X\rightarrow {X}^{{\,\prime}}$  is said to be {\it
equicontinuous at a point} $x_0 \in X$ if for every $\varepsilon >
0$ there is $\delta > 0$ such that ${d}^{{\,\prime}}
(f(x),f(x_0))<\varepsilon$ for all $f \in \frak{F}$ and $x \in X$
with $d(x,x_0)<\delta$. The family $\frak{F}$ is {\bf
equicontinuous} if $\frak{F}$ is equicontinuous at every point $x_0
\in X.$ It is known that every normal family $\frak{F}$ of mappings
$f:X\rightarrow {X}^{\,\prime}$ between metric spaces $(X,d)$ and
$\left({X}^{\,\prime},{d}^{\,\prime}\right)$ is equicontinuous. The
inverse conclusion is true whenever $(X,d)$ is separable and
$\left({X}^{\,\prime},{d}^{\,\prime}\right)$ is compact metric space
(see the version of the Arzela--Ascoli's theorem mentioned above in
\cite[section~20.4]{Va$_1$}).

\medskip
Let $D$ be a domain in ${\Bbb R}^n,$ $n\ge 2,$
$Q:D\rightarrow[0,\infty]$ be a Lebesgue measurable function. Denote
by $\frak{R}_{Q, \Delta}(x_0)$ the family of all ring
$Q$--homeomorphisms $f:D\rightarrow \overline{{\Bbb R}^n}$ at $x_0$
with $h(\overline{{\Bbb R}^n} \backslash f(D)) \ge \Delta
>0.$ Let $\frak{R}_{Q, \Delta}(D)$ denotes a family of all homeomorphisms
$f:D\rightarrow \overline{{\Bbb R}^n}$ such that $f\in\frak{R}_{Q,
\Delta}(x_0)$ at every point $x_0\in D.$ Let us consider
$\frak{R}_{Q, \Delta}(x_0)$ as the family of the mapping between
metric spaces $(X, d)$ and $\left(X^{\,\prime}\,,
d^{\prime}\right),$ where $X=D,$ $X^{\,\prime}=\overline{{\Bbb
R}^n},$ $d(x,y)=|x-y|$ be Euclidean metric and
$d^{\,\prime}(x,y)=h(x,y)$ be chordal metric. The following
statement take a place (see \cite[Theorem~6.1, Theorem~6.5,
Corollary~6.7]{RS}).

\begin{proposition}\label{pr8}{\sl
The family $\frak{R}_{Q, \Delta}(x_0)$ is equicontinuous at $x_0\in
D$ whenever at least one of the conditions holds: 1) The relation
$\int\limits_{0}^{\varepsilon(x_0)}\frac{dt}{tq_{x_0}^{1/(n-1)}(t)}=\infty$
take a place for some $\varepsilon(x_0)>0,$ $\varepsilon(x_0)<{\rm
dist\,}(x_0,
\partial D);$
2) The relation $q_{x_0}(r)\le C\cdot\log^{n-1}\frac{1}{r}$ holds as
$r\rightarrow 0$ and some $C>0;$ 3) $Q\in FMO(x_0).$ Besides of
that, if at least one of the conditions 1)--3) holds for every
$x_0\in D,$ the family $\frak{R}_{Q, \Delta}(D)$ is equicontinuous
(normal) in $D.$}
\end{proposition}

\medskip
Denote by $\frak{F}_{Q, \Delta}(x_0)$ the family of all local ring
$Q$--homeomorphisms $f:D\rightarrow \overline{{\Bbb R}^n}$ at $x_0$
with $h(\overline{{\Bbb R}^n} \backslash f(D)) \ge \Delta
>0.$ Let $\frak{F}_{Q,
\Delta}(D)$ denotes a family of all local homeomorphisms
$f:D\rightarrow \overline{{\Bbb R}^n}$ such that $f\in\frak{R}_{Q,
\Delta}(x_0)$ at every point $x_0\in D.$ The following statement
taka a place.

\medskip
\begin{theorem}\label{th3}
{\sl The family $\frak{F}_{Q, \Delta}(x_0)$ is equicontinuous at
$x_0\in D$ whenever at least one of the conditions holds: 1) The
relation
$\int\limits_{0}^{\varepsilon(x_0)}\frac{dt}{tq_{x_0}^{1/(n-1)}(t)}=\infty$
take a place for some $\varepsilon(x_0)>0,$ $\varepsilon(x_0)<{\rm
dist\,}(x_0,
\partial D);$
2) The relation $q_{x_0}(r)\le C\cdot\log^{n-1}\frac{1}{r}$ holds as
$r\rightarrow 0$ and some $C>0;$ 3) $Q\in FMO(x_0).$ Besides of
that, if at least one of the conditions 1)--3) holds for every
$x_0\in D,$ the family $\frak{F}_{Q, \Delta}(D)$ is equicontinuous
(normal) in $D.$}
\end{theorem}

\medskip
\begin{proof} It follows from the assumptions that every mapping
$f\in \frak{F}_{Q, \Delta}(x_0)$ omits at least two values $a_f$ and
$b_f$ in $\overline{{\Bbb R}^n}.$ Let $T_f$ be a M\"{o}bius
transformation mapping the $b_f$ to $\infty.$ Since the M\"{o}bius
transformations preserve the moduli of curve's families (see
\cite[Theorem~8.1]{Va$_1$}), the family of mappings
$\widetilde{\frak{F}}_{Q, \Delta}(x_0)=\{\widetilde{f}=T_f\circ f:
f\in \frak{F}_{Q, \Delta}(x_0)\}$ consists of the local ring
$Q$--homeomorphisms $f:D\rightarrow {\Bbb R}^n$ at $x_0,$ which omit
$\widetilde{a_f}\in {\Bbb R}^n$ and $\infty.$ Suppose that one of
the cases 1), 2) or 3) take a place. Then every $\widetilde{f}\in
\widetilde{\frak{F}}_{Q, \Delta}(x_0)$ is injective in some
neighborhood of $x_0$ whose radius depends only on $n$ and $Q$ (see
Theorem \ref{th1A}, Theorem \ref{th2} or Corollary \ref{cor1},
correspondingly). Now the family $\widetilde{\frak{F}}_{Q,
\Delta}(x_0)$ as well as the family $\frak{F}_{Q, \Delta}(x_0)$ is
equicontinuous at $x_0$ by Proposition \ref{pr8}. The corresponding
conclusion for the family $\frak{F}_{Q, \Delta}(D)$ follows from the
proved above.
\end{proof} $\Box$

\medskip
\begin{rem}\label{rem1}
The results of the paper are not true for $n=2$ that shows the
example $f_m(z)=e^{mz},$ $m\in {\Bbb N},$ $z\in {\Bbb B}^2.$
\end{rem}

%=================Список литературы====================
%\end{fulltext}

\large
{\bf \noindent Evgenii A. Sevost'yanov} \\
Institute of Applied Mathematics and Mechanics,\\
National Academy of Sciences of Ukraine, \\
74 Roze Luxemburg str., 83114 Donetsk, UKRAINE \\
Phone: +38 -- (062) -- 3110145, \\
Email: brusin2006@rambler.ru
\end{document}